\documentclass[12pt]{article}

\usepackage{paralist}
\usepackage{natbib}
\setcitestyle{comma,aysep={}}
\usepackage{epsf}
\usepackage[pdftex]{graphicx}
\usepackage{amssymb,amsmath}
\usepackage{lscape}
\usepackage{umoline}
\usepackage{color}
\usepackage{multirow}
\usepackage{booktabs}
\usepackage{tikz}
\usetikzlibrary{shapes,arrows}
\usepackage[section]{placeins}
\usepackage{graphicx}
\usepackage{caption}   
\usepackage{subcaption}
\usepackage{units}
\usepackage{bm}  %bold the symbol in formulation
\usepackage[vlined]{algorithm2e}
\usepackage{secdot}

\usepackage{longtable} %table for multi pages

\usepackage[a4paper]{geometry}   %control the margin of page

\usepackage{tikz}
\usetikzlibrary{positioning}
\usetikzlibrary{arrows}
\usetikzlibrary{shapes.arrows}
\usetikzlibrary{matrix}

\newcommand{\myfigure}[3]{ \begin{figure}[tbp]
    \begin{center}
    #1
    \end{center}
    \caption{#2} \label{#3} \vskip -0.14in \end{figure} }

\newcommand{\mytable}[3]{ \begin{table}[tbp]
  \begin{center}
  \caption{#2} \vspace{4mm}
    #1
  \label{#3} \vskip -0.01in \end{center} \end{table} }

\usepackage{changes}
\definechangesauthor[name={Weili Zhang}, color=red]{WZ}
\setremarkmarkup{(#2)}

\begin{document}

\baselineskip 20pt plus .3pt minus .1pt

\begin{center}

{\LARGE Objective Scaling Ensemble Approach for Integer Linear Programming }\\[12pt]

% Authors and addresses:

\footnotesize
\mbox{\large Weili Zhang$^{a}$, Charles D. Nicholson$^{a,*}$}\\
$^{a}$ School of Industrial and Systems Engineering, University of Oklahoma, Norman, OK, USA.\\
$^{*}$ Corresponding author: Charles Nicholson,  \mbox{cnicholson@ou.edu}.\\

\normalsize

\end{center}

\noindent Abstract: The objective scaling ensemble approach is a novel two-phase heuristic for integer linear programming problems shown to be effective on a wide variety of integer linear programming problems.  The technique identifies and aggregates multiple partial solutions to modify the problem formulation and significantly reduce the search space. An empirical analysis on publicly available benchmark problems demonstrate the efficacy of our approach by outperforming standard solution strategies implemented in modern optimization software.  

\bigskip

% Here are the Keywords:

% Here is the History:
\noindent {\it History:} Withdraw from European Journal of Operations Research on Oct 29, 2017. Submit to Journal of Heuristics on Dec 6, 2017.

\noindent\hrulefill
% The body of the paper starts here:

\newpage
\noindent {\LARGE Objective Scaling Ensemble Approach for Integer Linear Programming }\\[12pt]

\noindent Abstract: The objective scaling ensemble approach is a novel, two-phase heuristic for integer linear programming problems shown to be effective on a wide variety of integer linear programming problems.  The technique identifies and aggregates multiple partial solutions to modify the problem formulation and significantly reduce the search space. An empirical analysis on publicly available benchmark problems demonstrate the efficacy of our approach by outperforming standard solution strategies implemented in modern optimization software. 

\section{Introduction}
\label{sec-Introduction}

Integer programming (IP) is a fundamental approach to NP-hard combinatorial problems that arise in wide range of application areas including production, scheduling, finance, network design, and others \citep{nemhauser1988integer, zhang2017bridge, zhang2018probabilistic, zhang2010lattice, zhang2016multi, zhang2010reformed, zhang2017resilience, zhang2016resilience}. There are both linear and non-linear formulations of IP problems.  In this investigation we focus on the former.  Broadly defined, an mixed integer linear program (MILP) aims at optimizing a linear objective function (without loss of generality we assume minimization) subject to a set of linear equality/inequality constraints over real and integer/binary variables.  Borrowing notation from \citet{fischetti2003local}, we define the MILP problem as %\eqref{IPmin}-\eqref{IPvar},
%\begin{equation}
\begin{align} 
	\text{(ILP)} \quad  \text{min} \  z(\mathbf{x}) &= \mathbf{c}^T\mathbf{x} \label{ILP1}  \\
	  A\mathbf{x} \geq \mathbf{b}  & \label{ILP2}\\
	  x_{j} \ \text{integer}  & \quad \forall{j \in \mathcal{G}}  \label{ILP3}\\
	  x_{j} \in \left\lbrace 0,1 \right\rbrace & \quad \forall{j \in \mathcal{B}}  \label{ILP4} \\
	  x_{j} \ge 0 & \quad \forall{j \in \mathcal{N}} \label{ILP5}.
\end{align} 
%\end{equation}

\noindent Here, $\mathbf{c}$ is an $n$-dimensional vector of costs, $\mathbf{x}$ is an $n$-dimensional vector of decision variables, $A$ is an $m \times n$ constraint  matrix, $\mathbf{b}$ is an $m$-dimensional vector of parameters, and $\mathcal{N}$ is a set of variable indices $\{1 \ldots n \}$ partitioned into three sets, $\mathcal{N} = \{\mathcal{B},\mathcal{G},\mathcal{C}\}$ associated with binary, integer, and continuous variables, respectively.
%where  $\mathcal{B}$ associated with binary variables, $\mathcal{G}$ associated with integer variables, and $\mathcal{C}$ associated with continuous variables.  
\noindent If $\mathcal{C}$ and $\mathcal{I}=\{\mathcal{B}, \mathcal{G}\}$ are non-empty, the problem is a \emph{mixed integer} linear program (MILP). If $\mathcal{C} = \emptyset$ and $\mathcal{G} \neq \emptyset$, the problem is a \emph{pure integer} problem.  If $\mathcal{C}$ and $\mathcal{G}$ are empty, but $\mathcal{B} \neq \emptyset$, it is a \emph{binary} programming (BP) problem. If $\mathcal{G}$ is empty, but $\mathcal{C}$ and $\mathcal{B}$ are not, it is a \emph{mixed binary} problem (MBP).  Finally, if $\mathcal{I}$ is empty, the problem is not an integer problem but a linear programming (LP) problem.  Various IP solution approaches entail temporarily removing the integrality constraints in \eqref{ILP3} and \eqref{ILP4} solving the associated LP relaxation.

Commercial MILP solvers, such as CPLEX and Gurobi used in both academia and industry, leverage branch-and-bound and cutting planes algorithms with linear programming relaxation to find exact optimal solution \citep{lodi2010mixed}. Due to the time and/or resource complexity of finding exact MILP solutions, there is value in obtaining near-optimal solutions to such problems quickly.  A large body of research has been directed towards finding solution approaches applicable to particular subclasses of MILP problems, e.g. fixed-charge network flow problems \citep{bertsimas2003robust}, network design \citep{Crainic00}, vehicle routing \citep{Gulczynski2011794}, and scheduling \citep{hoffman1993solving, Akker2000, belien2007exact}.  General-purpose MILP solution approaches on the other hand include ``pivot and complement'' for BP problems \citep{balas1980pivot}, ``pivot and shift'' for MILP problems \citep{balas2004pivot}, ``pivot, cut and dive'' \citep{eckstein2007pivot}, OCTANE for BP problems \citep{balas2001octane}, relaxation induced neighborhood search (RINS) \citep{danna2005exploring}, local branching \citep{fischetti2003local}, feasibility pump \citep{fischetti2005feasibility,achterberg2007improving,bonami2009feasibility}, and others  \citep{blum2003metaheuristics, patel2007active}.  Approximate solutions may be of sufficient quality to stand on their own or be used in combination with an exact procedure to find feasible solutions.  Modern solvers incorporate many heuristics as  part of the overall optimization strategy to improve time to solution \citep{linderoth2010milp}.  

%\citep{land1960automatic,dakin1965tree}

In this paper, we introduce an heuristic technique and basic framework suitable for a wide variety of MILP problems.  This algorithm, which we call the  ensemble approach (OSEA), is inspired by RINS \citep{danna2005exploring} and slope scaling \citep{kim1999solution} techniques. Section \ref{OSEAF} reviews these techniques, explains the  motivation of OSEA, and formally defines the framework.  Section \ref{results} describes the problem testbed  and reports the computational results.  We summarize the work in Section \ref{conclusion}.

\section{Objective Scaling Ensemble Approach}
\label{OSEAF}

%In this section, we review the motivating literature and concepts for OSEA, discuss some preliminary analysis, and explain the framework in detail . Finally, two parameters will be introduced to tune OSEA performance with respect to balancing solution efficiency and  quality.

\subsection{Motivation} \label{motivation}

Relaxation induced neighborhood search \citep{danna2005exploring} is one of several heuristic techniques used in conjunction with exact solution approaches to MILP problems.  Branch-and-bound (or branch-and-cut) explore the MILP solution space by iteratively fixing one or more integer variables (a partial solution) and solving the remaining subproblem as a linear relaxation. The partial solutions are typically referred to as nodes in the search tree.  The best integer feasible solution found during the process is called an incumbent solution and is updated whenever a better feasible solution is found.  The process is iterated until the entire solution space has been implicitly examined and a provably optimal solution is found (assuming a feasible solution exists).  Many heuristic search techniques, on the other hand, search a \emph{neighborhood}, a local space ``close'' to a particular point within the solution space (as defined by some distance measure), to find improved solutions, e.g. local branching \citep{fischetti2003local}, tabu search \citep{glover1999tabu}, simulated annealing \citep{SA1983}, and machine learning \citep{zhang2016prediction,nicholson2016optimal} See \citet{gendreau2010handbook} for an excellent resource regarding a wide variety of metaheuristic techniques.  The local search is then repeated based on the neighborhood of the improved solution.  Neighborhood search procedures are often terminated based on some pre-specified criteria and do not guarantee the global optimality of the final solution. 

RINS employs information from the branch-and-bound process to form a search neighborhood of an incumbent feasible solution.  The intuition is that some subset of variables in a linear relaxation for a given search node will share values with the current incumbent solution.  The variables which do agree are fixed to their incumbent values.  The solution space for the resulting subproblem defines the neighborhood of the incumbent and this space is searched using an exact technique.  Any integer feasible solution found is a globally feasible solution and possibly will improve the incumbent solution. The sub-IP problem is potentially large and some stopping criterion is used to terminate the local search.  The master branch-and-bound process is resumed with a potentially improved incumbent.  RINS can be employed at any search node.  At each node the LP relaxation may result in a different solution and thus the overall search is diversified.

Slope-scaling, and in particular the dynamic slope scaling procedure (DSSP) \citep{kim1999solution}, was originally designed for the fixed-charge network flow (FCNF) problems and has been applied to a various problem types including the piecewise linear network flow problem \citep{kim2000dynamic}, the multicommodity fixed-charge network problem \citep{crainic2004slope}, the multicommodity location problem \citep{gendron2003tabu}, the minimum toll booth problem \citep{Bai2010}, and stochastic integer programming \citep{shiina2012DSSP}. With respect to the original application, DSSP employs a series of linearizations of the FCNF discontinuous objective function in \eqref{FCNFobj},
\begin{align} 
\hspace{1em}  \text{min} \hspace{0.5em}  & \sum_{(i,j) \in A} (c_{ij}x_{ij} + f_{ij}y_{ij})  \label{FCNFobj}
\end{align}
\noindent where $A$ is the set of arcs in the network, $\mathbf{x}$ is a vector of continuous arc flow values, $\mathbf{y}$ is a vector of binary decision variables, and the cost vectors $\mathbf{c}$ and $\mathbf{f}$ are the unit flow costs and fixed costs, respectively.  DSSP removes the binary variables $\mathbf{y}$ from the objective, iteratively scales the fixed cost $\mathbf{f}$ and adds them to the variable flow cost $\mathbf{c}$. Similar to RINS, OSEA uses information from linearized formulations to form a smaller search space.  The relaxation formulation differs from RINS in that the relaxation is not based on a search node (partial solution), but based on a series of slope-scaling inspired relaxations.  

OSEA defines the search space based on  variable agreement between (i) an incumbent solution and one or more solutions of the linearized formulations (i.e., an ensemble of solutions) or (ii) entirely from agreement between solutions in the ensemble.  Solution ensembles have been exploited in a variety of ways such as variable fixing based on value agreement among solutions (e.g., as in RINS) or  ``voting''  among the solutions (e.g., 4 of 5 solutions have variable $x_1 = 17$).  The latter is commonly employed in the field of statistical learning \citep[e.g.,][]{breimanBagging96}.  OSEA takes a relatively conservative approach in using the ensemble to define a sub-MIP problem to be solved exactly.

\subsection{OSEA Framework}
In many instance of large, real-world IP problems, only a small percentage of the integer variables have non-zero values in the optimal solution. It is worth noting that among the larger MILP instances available in the IP benchmark problem library MIPLIB 2010 \citep{koch2011miplib}, the relative number of non-zero integer variables in the optimal solutions is very low.  Table \ref{MIPLIB1} shows the summary statistics (minimum, first quartile, median, mean, third quartile, and maximum) for the percentage of integer variables used in the optimal solution for the 36 MIPLIB 2010 solved problems that have at least 10,000 integer variables.  The median is 1.3\% and 75\% of the optimal solutions use less than 2.89\% of the possible integer variables.  Motivated by this feature of IP problems, OSEA attempts to \emph{eliminate} integer variables from the problem formulation and the ensemble aggregation method for OSEA is designed with this characteristic in mind. 

\mytable{
\begin{tabular}{c c c c c c }
\toprule
   Min & Q1 & Median &  Mean & Q3 &  Max \\ 
 0.01\%  & 0.28\% & 1.30\% &  3.69\% &  2.89\% & 34.41\% \\ 
\bottomrule
\end{tabular}}
{Integer Variables Used in Large MIPLIB instances}{MIPLIB1}

Let $\mathcal{E} = \left\{ \mathbf{s}_1, \mathbf{s}_2, \ldots, \mathbf{s}_k \right\}$ denote an ensemble of $k$ solutions (possibly including infeasible solutions) to an MILP problem.  The set of solutions may be generated through slope-scaling techniques, LP relaxations, known feasible solutions, accumulating incumbent solutions in a branch-and-bound algorithm, or other methods.  In particular, OSEA fixes the $j^\text{th}$ integer variable to 0 if for all solutions $s \in \mathcal{E}$, the $j^\text{th}$ variable, $s_j$, equals 0.  That is, 
\begin{align*}
x_j \leftarrow 
\begin{cases} 
 \text{fix to 0} \quad & \text{if } s_j = 0 \quad \forall s \in \mathcal{E} \\
 \text{do not fix}  & \text{otherwise}
 \end{cases} \quad \forall j \in \mathcal{I}. 
\end{align*} 
\noindent The integer variables which are left open in the corresponding sub-MIP problem form a reduced search space.  An exact search of the reduced problem space produces the OSEA solution and objective.

It is important to note that the ensemble does not necessarily consist of high quality solutions to the original problem.  In fact, from initial testing we place a priority on diversity of quality.  If the solutions in $\mathcal{E}$ are diverse, then the variables that are unused by every solution in the ensemble share at least one characteristic: they are each ``unattractive'' to a wide range of solutions.  Since OSEA fixes variables to 0, by allowing poorer quality solutions in the ensemble, we take a more conservative approach.  That is, only variables that are not used among a variety of solutions (e.g., good, median, poor) are discarded.

Moreover, if a given integer variable is not used in a linearized optimal solution when the cost is adjusted to a fraction of its original cost, then the intuition is that the integer variable is not likely useful in the original problem. The absolute cost associated with a variable is not as important as the cost relative to other variables.  The iterative scheme in slope-scaling techniques and in OSEA updates individual variable costs throughout the process.  This update scheme dynamically affects the relative costs of the variables.  Variables which may be too costly to use in a linearized solution during the earlier iterations may become cost-effective in the latter ones.

The ensemble $\mathcal{E}$ must be populated with solutions or partial solutions.  For OSEA, this is primarily accomplished in the objective scaling iteration phase. The scaling process is now described.

OSEA scales the coefficients of the discrete variables and iteratively solves the relaxed problem \eqref{objLP}, 

\begin{equation}
\begin{aligned} 
	(\text{LP}_n)  \ \ & \text{min} \ z^\text{LP}_n(\mathbf{x}) = \sum_{j \in \mathcal{C}}c_{j}x_{j} + \sum_{j \in \mathcal{I}} \bar{c}^{n}_{j}x_{j}  \\
	   A\mathbf{x} & \geq \mathbf{b}   \\
	   0 \leq x_{j} &\leq 1  \quad \forall{j \in \mathcal{B}}   \\
	   x_{j} & \ge 0  \quad \forall{j \in \mathcal{N}}
\end{aligned} \label{objLP}
\end{equation}

\noindent where $\bar{c}^{n}_j$ for $j \in \mathcal{I}$ is the scaled cost coefficient of the integer variable. Let $\tilde{\mathbf{x}}^{n}$ denote the solution to $\text{LP}_n$.  This solution is used to update the integer coefficient for the next iteration $n+1$ as follows,

\begin{align}
\bar{c}_{j}^{n+1} \leftarrow \frac{c_{j}}{\tilde{x}_{j}^n+1} \quad \forall j \in \{\mathcal{I}:\tilde{x}_{j}^n > 0 \}  \label{eq4}
\end{align}

Note in DSSP, the fixed cost value is scaled by $\nicefrac{1}{\tilde{x}_{j}^n}$, however there are two benefits that result from modifying this for general application.  First, since MILP problems may have negative integer variable cost coefficients, as the relaxed solution approaches 0, the scaled costs may approach negative infinity, $$\lim_{\tilde{x}^n_j \to 0} \frac{-1}{\tilde{x}_{j}^n} = -\infty $$and a counterintuitive result ensues, namely the attractiveness of the variable increases without bound as the value of the variable decreases.  This effect is bounded by a simple modification of the denominator in Equation \ref{eq4}. And secondly, the resulting bound is intuitive 
$$\lim_{\tilde{x}^n_j \to 0} \frac{c_j}{\tilde{x}_{j}^n + 1} = c_j.$$ 
 
Slope-scaling dynamically modifies the costs of different variables throughout the search process to alter their relative attractiveness in the relaxations.  While it is true that modification in Eq. \eqref{eq4} impacts the appealing characteristic of the final iteration $N$ of DSSP in which the scaled objective value reflects the true solution cost (i.e., it includes the full fixed cost incurred in the corresponding network flow solution),
$$
\sum_{j \in \mathcal{C}}c_{j}\tilde{x}^N_{j} + \sum_{j \in \mathcal{I}} \bar{c}^{N}_{j} [\tilde{x}^N_j > 0]
$$
\noindent where $[x^N_j > 0]$ denotes the Iverson bracket which returns a 1 if $x^N_j > 0$ and 0, otherwise.  However, this outcome is not critical to the success of DSSP.  That is, the best solutions from DSSP are not necessarily found in the final iteration \citep{Nahapetyan2008}.   At earlier iterations $n < N$, the scaled fixed costs do not represent the true value: $$\sum_{j \in \mathcal{I}} \bar{c}^{n}_{j}\tilde{x}^n_{j}  \ne \sum_{j \in \mathcal{I}} \bar{c}^{N}_{j} [\tilde{x}^N_j > 0].$$ This suggests that the search path induced by the procedure is of more importance than the objective value of the final iteration.

The integer variable cost coefficients are initialized to a fraction of the original cost by scaling by the inverse of the relatively large value $M$, 
\begin{align*}
 M = \sum_{j \in \mathcal{I}}|c_{j}|.
\end{align*} 
\noindent Thus,  
\begin{align*}
 	c_{j}^{0} \leftarrow \frac{c_{j}}{M} \qquad  \forall j \in \mathcal{I}
 \end{align*}
and the update scheme is then,
\begin{align*}
 	 c_{j}^{n} = \begin{cases}
 	 \frac{c_{j}}{\tilde{x}_{j}^{n-1}+1} \qquad & \text{if } \tilde{x}_{j}^{n-1} > 0\\[2ex]
 	  c_{j}^{n-1} \qquad & \text{otherwise} \\
 	 \end{cases} \qquad \forall j \in \mathcal{I}.
\end{align*} 

In the objective scaling phase of OSEA there are $N>0$ iterations and consequently $N$ linearized solutions: a subset of which will be added to the ensemble $\mathcal{E}$. For large $N$, the number of different solutions that could be added to $\mathcal{E}$ is also large.  If $|\mathcal{E}|$ is too large, the reduced search space is potentially too large for practical purposes.  Therefore, we select a subset of the iterated linear solutions to be added to the ensemble.  Let $\mathcal{S}$ denote this subset.  On the other hand, if $|\mathcal{E}|$ is too small or if it does not contain sufficient diversity, then there may be insufficient options in the search space to generate good solutions to the original MILP.  A number of possible strategies can be designed to build the set $\mathcal{S}$.  We devise one such possible strategy to emphasize ensemble diversity in Section \ref{results}.

If the ensemble $\mathcal{E}$ does not contain a feasible solution then OSEA may or may not produce a feasible solution.  However, if a feasible solution is included in $\mathcal{E}$ then OSEA is guaranteed to find a feasible solution in the reduced search space. In our implementation of the framework, we take a hybrid approach in which we utilize the already existing pre-processing, heuristics, and branch-and-cut algorithms readily available in commercial software such as Gurobi and CPLEX to briefly search for a feasible solution that can be added to $\mathcal{E}$. This will be described in more detail in Section \ref{results}.  

There are multiple possible stopping criterion for OSEA.  Similar to DSSP, the iterative objective scaling phase will stop once there are no new cost coefficient updates for the integer variables. In some cases it might be prudent to provide an upper limit on the total number of iterations allowed for the iterative scaling procedure.  Let $N_\text{max}$ denote the max allotted iterations.  Additionally, since OSEA is meant as an heuristic technique to quickly reduce the search space of complex problems, a time limit could also be imposed on the scaling phase. The complete OSEA logic (including the iteration limit, but not the time limit) is summarized in Figure \ref{OSEACODE}.

%\mytable{
%\begin{tabular}{c c}
%$\beta$ & Iterations \\
%\toprule
%1.0 & 221 \\
%1.5 & 140 \\
%2.0 & 60 \\
%\bottomrule
%\end{tabular}}
%{Iterations till convergence for \emph{dc1c} FIX THIS!(Updated)}{beta}
\myfigure{
	\begin{algorithm}[H]
		\DontPrintSemicolon
		\KwData{ILP instance, $N_\text{max}$}
		\KwResult{$z_{\text{OSEA}}$, $\mathbf{x}_\text{OSEA}$}
		\Begin{
		    compute $M$:   $M \gets \sum_{j \in \mathcal{I}}|c_{j}|$ \;
		    initialize integer costs : 	$c_{j}^{0} \gets \frac{c_{j}}{M} \  \forall j \in \mathcal{I}$\;
		    optional: initialize $\mathcal{E}$ by including one or more feasible solution(s) to MILP \;
					$n \gets 0$ \;
					\While{$n \leq N_\text{max}$}
					{   \vspace{1mm}
						$\mathbf{\tilde{x}}^n \gets $ solution to Problem $\text{LP}_n$ \;
						\vspace{1mm} $\bar{c}_{j}^{n} \gets \begin{cases}
						 	  \displaystyle\frac{c_{j}}{\left(\tilde{x}_{j}^{n-1} +1\right)}\ & \text{if } \tilde{x}_{j}^{n-1} > 0\\[2ex]
						 	  \bar{c}_{j}^{n-1} \ & \text{otherwise} \\
						 	 \end{cases} \ \forall j \in \mathcal{I}$ 
						 \vspace{1mm} \;
						 	
						 \If{$\bar{c}_{j}^{n} = \bar{c}_{j}^{n-1} \ \forall j \in \mathcal{I}$}  
						 { \vspace{1mm}
						 	\textbf{break}
						 }				 
					 
					$n \gets n +1$ \;
					}
					 $\mathcal{S} \gets $ a subset of $\{ \mathbf{\tilde{x}}^1, \dots, \mathbf{\tilde{x}}^N \}$ \;
					  
					 $\mathcal{E} \gets \mathcal{E} \cup \mathcal{S}$ \;
				
		$x_j \leftarrow 
		\begin{cases} 
		 \text{fix to 0} \quad & \text{if } s_{j} = 0 \ \forall s \in \mathcal{E} \\
		 \text{do not fix}  & \text{otherwise}
		 \end{cases} \quad \forall j \in \mathcal{I} $ \;
		$z_{\text{OSEA}}$, $\mathbf{x}_\text{OSEA} \gets$ solve reduced MILP problem \;
		}
	\end{algorithm}
}{Objective Scaling Ensemble Approach for MILP problems}{OSEACODE}

\section{Computation Results}
\label{results}

\subsection{Experimental Design}
\label{testbed}

The MIPLIB 2010 \citep{koch2011miplib} is a publicly available library of  pure and mixed integer programming problem instances assembled by researchers and practitioners over several years.  This library of benchmark problems is used in evaluating software performance of commercial solvers \citep{gurobiMIPLIB2010}.  The library contains 361 instances classified into $3$ difficulty levels: $185$ \emph{easy}, $42$ \emph{hard}, and $134$ \emph{open} problems.  The latter problem class contains the instances which have yet to be solved optimally. The instances are further described by $8$ characterizations types: \emph{benchmarks} (B) are solvable within 2 hours on a PC, \emph{infeasible} (I), \emph{primal} (P) instances have the LP relaxation objective equal to the optimal objective, \emph{extra-large} problems (X), \emph{reoptimize} (R) instances require a relatively long time to solve the LP relaxations, \emph{tree} (T) instances have a large number of enumeration trees, \emph{unstable} (U) instances have poor numerical properties, and \emph{challenge} (C) instances which are classified generally as difficult to solve. The majority of the instances in the library also include information relating to the problem application area (e.g. lot sizing, open pit mining, network design, etc.)  

OSEA is tested on a 170 problem subset of the 361 MIPLIB problems.  Since OSEA is appropriate only for problems with integer variables in the objective function, any problems without this characteristic are discarded (87 problems).  Infeasible problems are also eliminated from testing (22 problems).  Instances which exceed the memory capacity of our available computer equipment (23 problems) and those in which a feasible solution was not found within 60 seconds are not included in the experimentation (59 problems).

The experiments will be conducted as follows.  To increase the likelihood that an integer feasible solution is included in the ensemble $\mathcal{E}$, we will employ the existing pre-processing and heuristic algorithms of the optimization software by attempting to solve each problem for one second using commercially available optimization software.  If a feasible solution is found in the time limit, it will be added to the ensemble. Regardless, the scaling phase begins and we select solutions discovered during this iterative stage to be added to the ensemble. Based on initial testing we emphasize ensemble diversity by including three solutions from the scaling iterations into the ensemble: the solutions associated with the best, worst, and median MILP objective values.

The maximum running time for OSEA including initial one second search, objective scaling phase, and solving the reduced problems is set to 60 seconds.  To compare the OSEA solution quality with the state-of-the-art exact techniques, we will use the Gurobi optimization software version 5.6.3 with a time limit of 60 seconds.  According to the Gurobi product website, this solver includes 14 different MIP heuristics, 16 cutting plane strategies, and a variety of presolve techniques \citep{gurobiProduct}.  We use the default parameter settings for Gurobi with all heuristics activated (including RINS). The best objective value found using the default settings of the Gurobi optimization software and the time to find that value are recorded.  All tests are performed on a Windows 7 64bit machine with Intel Xeon CPU E5-1620 and 8 GB RAM with a single thread.  

OSEA is inherently dependent on an IP solver.  The technique iteself is used in conjuction with a solver to reduce the IP search space.  The empirical analysis will compare the results of using OSEA with a commercial solver against using the same commercial solver without OSEA.  While the commericial solver is not being tuned specifically for each problem, the settings are identical for the OSEA test.  That is, the only difference is the additional OSEA overhead to the commercial solver.

\subsection{Experimental Results}
\label{eresults}

Let $t_\text{OSEA}^{p}$ and $t_\text{standard}^{p}$ denote the computing time to solve problem $p\in P$ using OSEA and the standard (solver without OSEA), respectively.  Similarly, the best objective values found for problem $p \in P$ denoted by $z_\text{OSEA}^{p}$ and $z_\text{standard}^{p}$, respectively. Additionally, the solution gap (ILPgap) is used to evaluate solution quality.  Let $G_{X}^{p}$ denote the MILP gap for approach $X \in \{\text{OSEA}, \text{standard}\}$ on problem $p \in P$,

 $$G_{X}^{p}=\dfrac{|z_\text{bound}^{p}-z_\text{X}^{p}|}{|z_\text{X}^{p}|} \times 100\%, \ \forall X \in \{\text{OSEA, standard}\},$$

\noindent where $z_{bound}^{p}$ is the known optimal objective value for \emph{Easy} and \emph{Hard} problems, and is the linear relaxed objective value for \emph{Open} instances.  Note if $z_\text{X}^{p} = 0$, a small positive value is added to the denominator. 
 %\add{The experimental results are reported in Appendix \ref{Easy}, \ref{Hard}, and \ref{Open} based on the problem difficulty.}

Let $\gamma$ denote the percentage of integer variables removed in the objective scaling phase of OSEA. The distribution of $\gamma$ depicted in Figure \ref{distgamma} shows that OSEA removes a significant percentage of integer variables for the majority of test bed.  Overall, OSEA removes an average of 61.97\% of the integer variables. Reducing the MILP solution space can lead to notable improvements in computation time.

The results are summarized in Table \ref{averagePerformance} where $n$ equals the number of problem instances in each cell.  Overall, using a paired $t$-test, the difference in computing time and objective values of OSEA against the baseline are statistically significant.  Both techniques find optimal solutions for 40\% of the \emph{Easy} and \emph{Hard} problems,  albeit the instances differ. That is, OSEA finds optimal solutions to certain instances that the standard approach failed to find within the time limit, and visa versa.  Among the \emph{Easy} and \emph{Hard} problems, OSEA terminates faster than the standard technique, and while the average optimality gap is smaller for OSEA, the solution quality differences are not significant at a 95\% confidence level.   For the \emph{Open} problems, the results are statistically significant: OSEA produces a higher quality solution within the time limit and does so in less time.

\myfigure{
		\includegraphics[width=.75\textwidth]{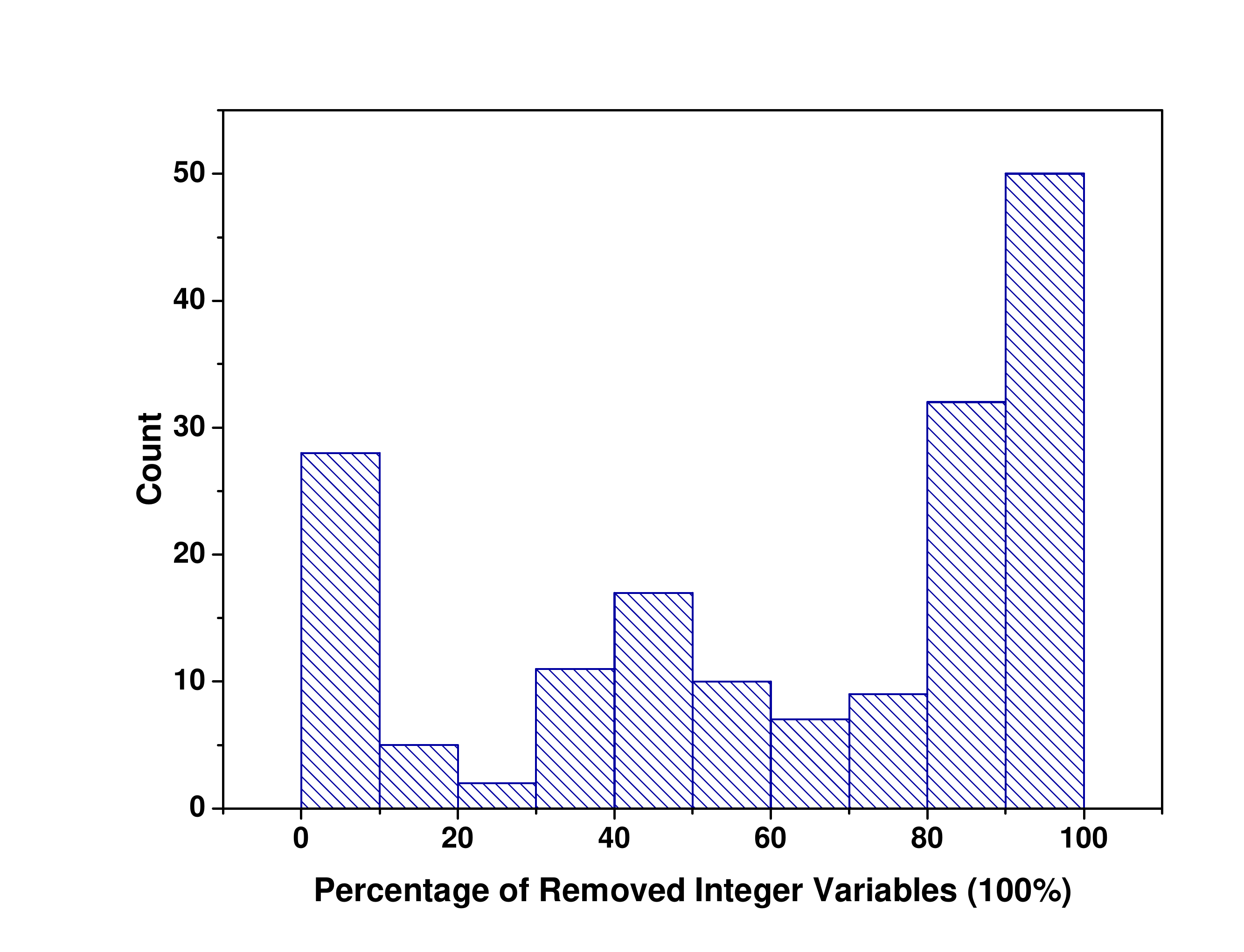}}
{Distribution of $\gamma$}{distgamma}

\mytable
{\begin{tabular}{l r r r r r r r }
\toprule
& & \multicolumn{3}{c}{Computing Time} & \multicolumn{3}{c}{Objective Value} \\
\cmidrule(lr){3-5} \cmidrule(lr){6-8}
Difficulty & $n$ & $t_\text{OSEA}$ & $t_\text{standard}$ &  $p$-value  & $G_\text{OSEA}$ &  $G_\text{standard}$ & $p$-value  \\
\toprule 
Easy & 108 & 19.84 & 41.41 & $<0.0001$ & 3.12\% & 4.9\% & 0.066 \\
Hard & 32 & 30.39 & 60.03 & $<0.0001$ & 8.4\% & 12.51\% & 0.072 \\
Open & 30 & 42.31 & 60.09 & 0.0004 & 1226.25\% & 1227.02\% & 0.028 \\
\midrule
Overall & 170 & 25.26 & 47.99 &  $<0.0001$ & 212.73\% & 214.77\% & 0.014 \\
\bottomrule
\end{tabular} }{Statistical Analysis of Performance by Difficulty Level}{averagePerformance}

We use the performance profile technique from \citet{dolan2002benchmarking} to further evaluate OSEA.  The baselines for comparisons on problem $p\in P$ are set as the best MILP gap and computing time, respectively. Let $r^t_{p,X}$ denote the performance ratio of computing time on problem $p\in P$ using technique $X$, $$r^t_{p,X} = \dfrac{t_{X}^{p}}{\text{min}\{t_{\text{OSEA}}^{p},t_{\text{standard}}^{p}\}}, \ \forall X \in \{\text{OSEA},\text{standard}\}.$$

Let $\rho_{X}^{t}(\tau)$ denote the probability for approach $X$ that $r^t_{p,X}$ is within a factor $\tau$ of the best ratio in terms of computing time, $$\rho^{t}_{X}(\tau)=\dfrac{1}{|P|}\text{size}\big\{p\in P: r^{t}_{p,X} \leq \tau \big\}, \ \forall  X \in \{\text{OSEA},\text{standard}\}.$$ Similarly, let  $\rho^{G}_{X}(\tau)$ denote the probability that approach $X$ is within a factor of $\tau$ from the best MILP gap ratio.  In both cases, larger values are preferred.  The cumulative distributions of $\rho^{t}_{X}(\tau)$ and $\rho^{G}_{X}(\tau)$ form the respective performance profiles.  Performance profiles evaluate the overall performance of a solution technique and when $|P|$ is sufficiently large, are relatively robust with respect to performance outliers of individual problem instances \citep{dolan2002benchmarking}.

\begin{figure}
 \centering
  \includegraphics[width=0.75\linewidth]{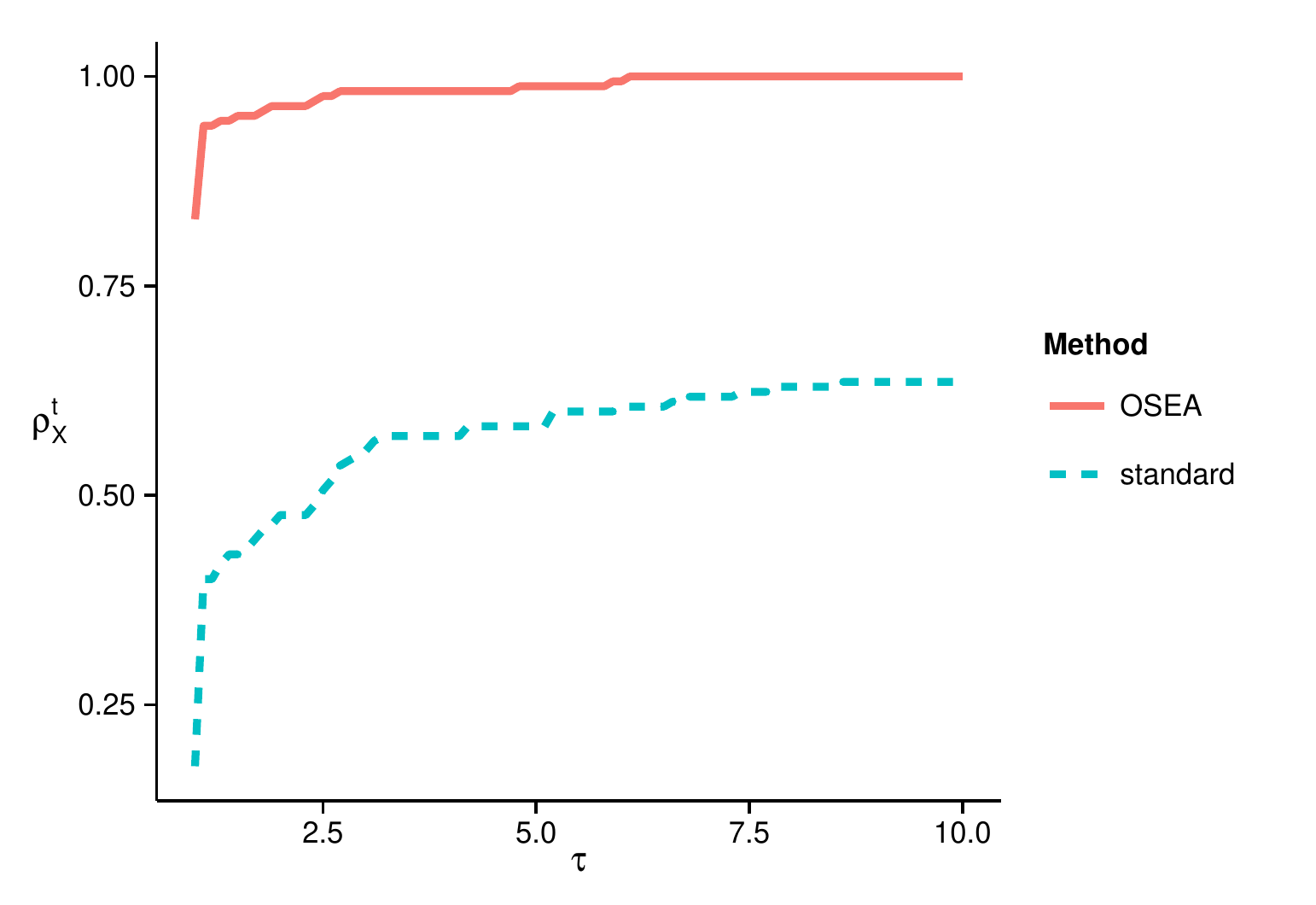}
  \caption{Computation time performance profile for $1 \le \tau \le 10$}
  \label{PPT1}
\end{figure}

\begin{figure}
  \centering
  \includegraphics[width=.75\linewidth]{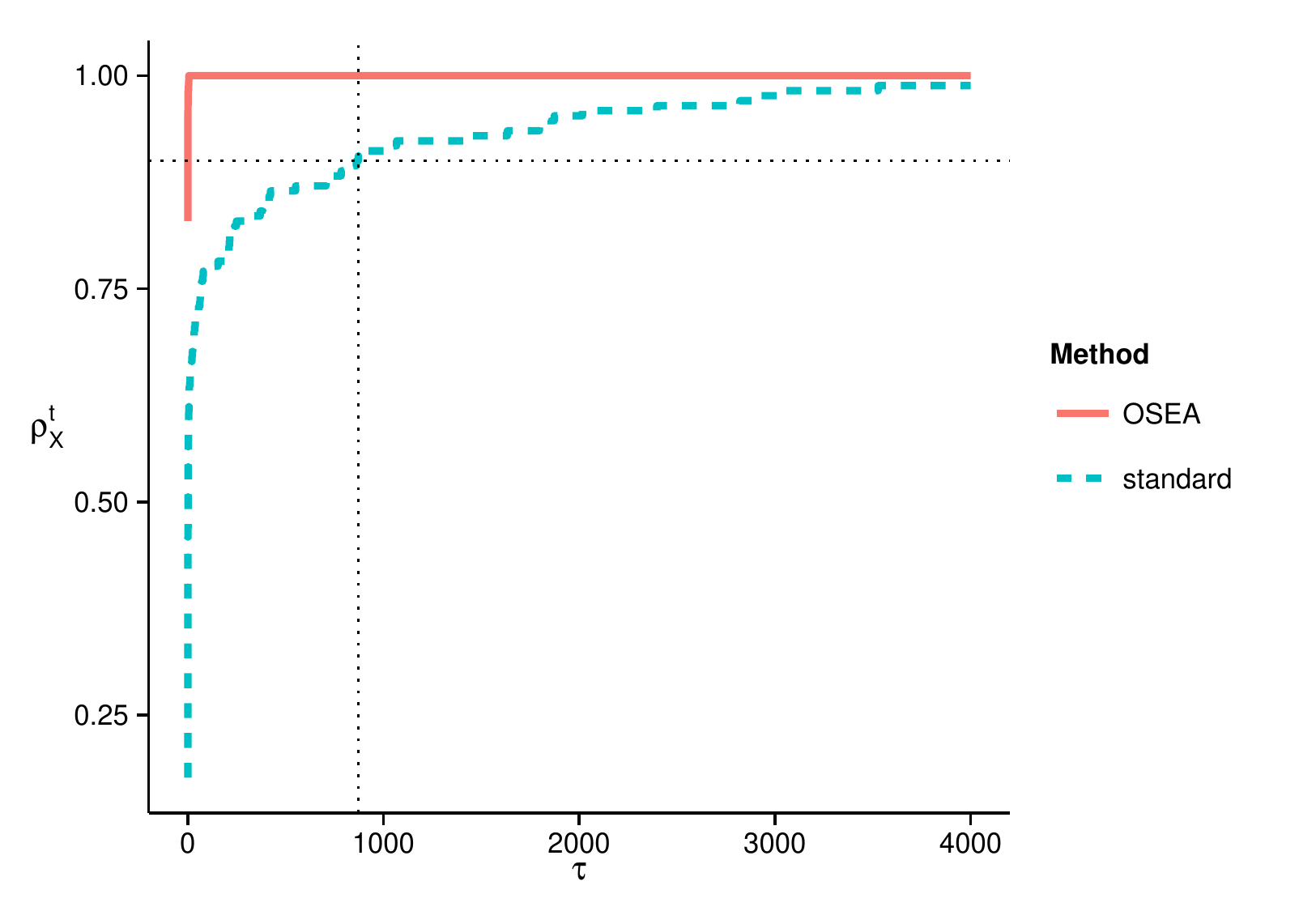}
  \caption{Computation time performance profile for $1 \le \tau \le 4000$}
  \label{PPT2}
\end{figure}

The performance profile for solution times are presented in Figures \ref{PPT1} and \ref{PPT2} with two different ranges for $\tau$ each.  The probability that OSEA terminates earlier than the standard approach is 0.829 (see Figure \ref{PPT1} when $\tau=1$).  OSEA solves 100\% of the problems within a factor of 6.1 for the computation time ratio, i.e.,  $\rho^{t}_{\text{OSEA}}(6.1)=1$, whereas  $\rho^{t}_{\text{standard}}(6.1)$ is only 0.606. The performance profile for the standard approach demonstrates that OSEA is much faster for many problems, e.g., $\rho^{t}_{\text{standard}}(\tau) \le 0.9$ for $\tau\le 872$ (depicted in Figure \ref{PPT2}.)  That is, OSEA terminates 872 times faster than the standard approach on 10\% of benchmark problems.  

\begin{figure}
  \centering
  \includegraphics[width=.75\linewidth]{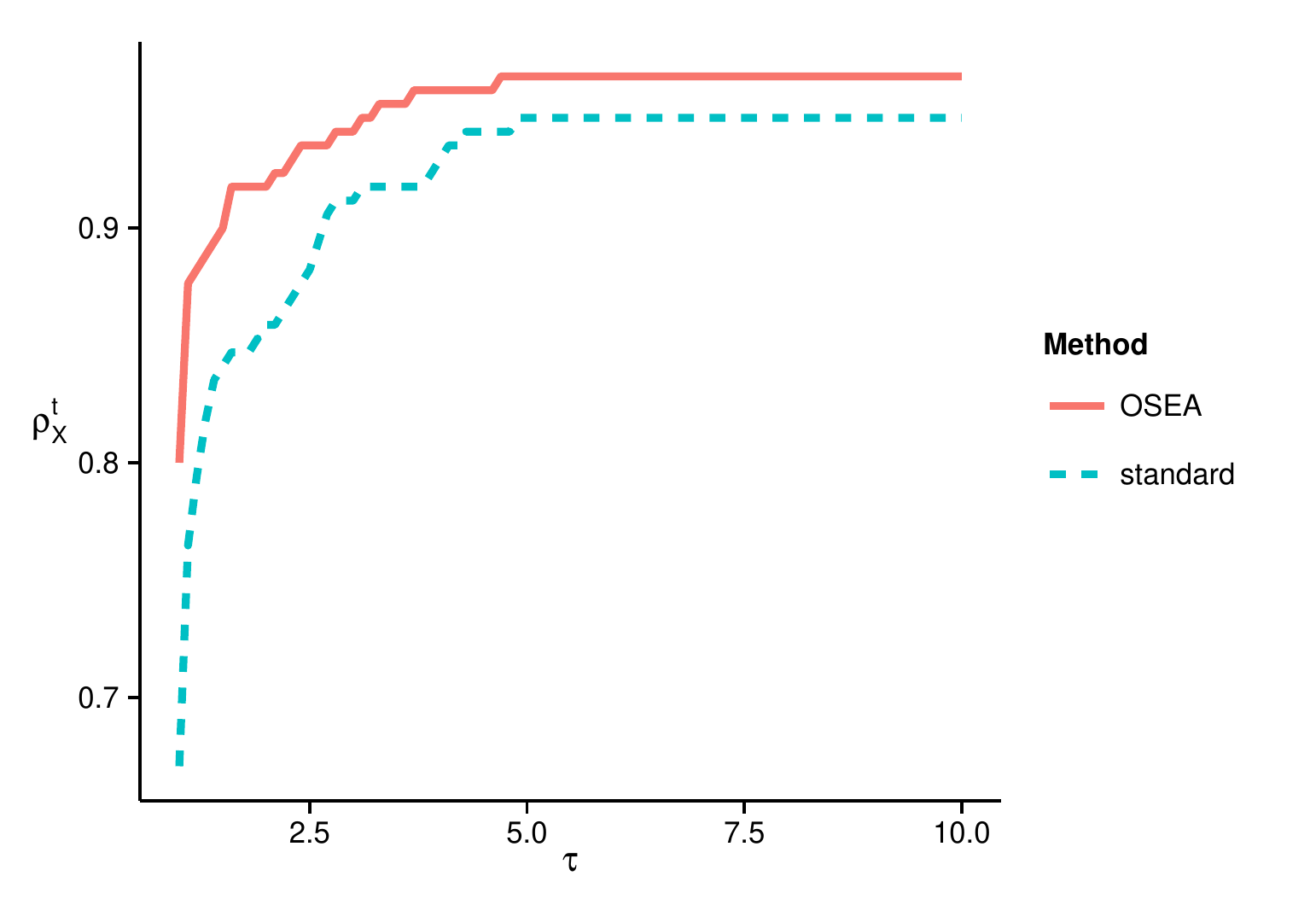}
  \caption{Solution quality performance profile for $1 \le \tau \le 10$}
  \label{PPG}
\end{figure}

While the reduced problems solve faster than the original problems, the solution qualities must be examined. Figure \ref{PPG} shows the performance profile of solution quality for $\tau \in [1,10]$.  OSEA is more likely to outperform the standard method, i.e., $\rho_\text{OSEA}^{G}(1) = 80\% > \rho_\text{standard}^{G}(1) = 67\%$. OSEA solves 90\% of all problems within a factor of 1.5 of the best technique. The OSEA solution quality performance profile is equal to or superior than the standard performance profile across all values of $\tau$. Note the minimum values of $\tau$ necessary to capture all problems: $\rho_\text{OSEA}^{G}(86) = 100\%$, whereas $\rho_\text{standard}^{G}(1173) = 100\%$ (not shown in figure) implying that OSEA performs, at worst, 86 times as bad as the standard technique, whereas the standard technique performed, up to 1173 times  worse than OSEA.

%\begin{figure}
%\begin{subfigure}{.5\textwidth}
%  \centering
%  \includegraphics[width=1.1\linewidth]{10pptime.pdf}
%  \caption{Performance profile on $[0,10]$}
%  \label{PPT1}
%\end{subfigure}%
%\begin{subfigure}{.5\textwidth}
%  \centering
%  \includegraphics[width=1.1\linewidth]{fullpptime.pdf}
%  \caption{Performance profile for All}
%  \label{PPT2}
%\end{subfigure}
%\caption{\added[id=WZ,remark=add]{Performance Profile on Computing Time}}
%\label{PPT}
%\end{figure}

%\begin{figure}
%\begin{subfigure}{.5\textwidth}
%  \centering
%  \includegraphics[width=1.1\linewidth]{10ppgap.pdf}
%  \caption{Performance profile on $[0,10]$}
%  \label{PPG1}
%\end{subfigure}%
%\begin{subfigure}{.5\textwidth}
%  \centering
%  \includegraphics[width=1.1\linewidth]{fullppgap.pdf}
%  \caption{Performance profile for All}
%  \label{PPG2}
%\end{subfigure}
%\caption{\added[id=WZ,remark=add]{Performance Profile on Solution Quality}}
%\label{PPG}
%\end{figure}

\section{Conclusions}
\label{conclusion}

The objective scaling ensemble approach is a novel, two-phase heuristic solution procedure that iteratively solves scaled linear versions of the original MILP problem and uses a subset of the LP relaxation results to form an ensemble of solutions.  This ensemble is aggregated in such a way to identify integer variables which are not likely to be used in an optimal solution.  These variables are removed from the MILP to create a reduced problem space.  Exact techniques such as branch-and-cut are applied to the revised problem formulation.  If the reduced search space is sufficiently large, a feasible and even possibly optimal solution for the original MILP can be found.  If the space is small enough, the revised problem space can be searched more efficiently.

The inspiration for OSEA comes from well known and successful heuristic approaches which have been used in conjunction with other techniques to produce a more efficient search of complex problem spaces.  Many advanced heuristic approaches are often invoked by default in commercial optimization software.  We compare the solution quality of OSEA in the first 60 seconds of optimization time to that of the assortment of heuristics, cutting plane strategies, and exact search algorithms implemented in Gurobi 5.6.3. OSEA successfully reduces the search space in a way which is competitive with industry leading optimization software.  

The empirical results on 170 publicly available benchmark integer programming problems and rigorous analysis indicate that OSEA can improve MILP solution quality on a wide range of problems without compromising the computation time.  Among the benchmark problems, many are well documented and related to published work \citep[e.g.,][]{fischetti2005feasibility, Bley20101641, RaackKosterOrlowskiWessaely2011}.  The instances include a wide variety of problem types and application areas including network design, open pit mine production, the p-Median problem, crew scheduling, and lot sizing, among others.  The problems range in size from hundreds of integer variables to several orders of magnitude more.  For certain problem types evaluated, OSEA performs exceptionally well, e.g. open pit mining.  For others, the results while promising, are mixed, e.g. network design problems.  In future work we will examine the particular nature of certain problem formulations to understand whether or not the outstanding results are generalizable to the problem class.

OSEA can be easily applied to any MILP problem with integer variables in the objective function.  However, OSEA is not meant to be used exclusive of traditional IP solvers, but ideally to be incorporated as yet another of the integrated heuristics used in software.   In our initial experimentation to this end, we find OSEA to improve the solution performance at the root node of the branch-and-bound algorithm, but not to work well at subsequent nodes. While some heuristics (e.g., RINS) are activated at various nodes of  branch-and-bound, earlier indications are that OSEA is a beneficial initial heuristic applied specifically at the root node to find better incumbent solutions early on. 

\newpage
\bibliographystyle{ijocv081}
\bibliography{references}

\color{blue}

\end{document}